\documentclass[letterpaper, 10 pt, conference]{ieeeconf}

\IEEEoverridecommandlockouts
\title{\LARGE \bf
	PIETOOLS: A Matlab Toolbox for Manipulation and Optimization of Partial Integral
	Operators
}

\author{Sachin Shivakumar$^{1}$, Amritam Das$^{2}$  and Matthew M. Peet$^{1}$
	\thanks{$^{1}$ Sachin Shivakumar\{sshivak8@asu.edu\} and Matthew M. Peet\{mpeet@asu.edu\} are with School for Engineering of Matter, Transport and Energy, Arizona State University, Tempe, AZ, 85298 USA }%
	\thanks{$^{2}$ Amritam Das\{am.das@tue.nl\} is with Department of Electrical Engineering, Eindhoven University of Technology}%
}

\let\labelindent\relax
\usepackage{enumitem}
\usepackage{graphicx}
\usepackage{amsfonts}
\usepackage{amssymb}
\usepackage{amsmath}
\usepackage{subcaption}
\usepackage{grffile}
\usepackage{multicol}
\usepackage{algorithm,algorithmic}
\usepackage{empheq}
\usepackage{mathtools}
\usepackage{xcolor}

\setlist[enumerate]{wide=0pt, widest=99,leftmargin=\parindent, labelsep=*}

\newcommand{\mcl}[1]{\mathcal{ #1}}
\newcommand{\mbf}[1]{\mathbf{ #1}}
\newcommand{\norm}[1]{\Vert #1\Vert}

\newcommand{\hinf}{\ensuremath{H_{\infty}}}
\newcommand{\ip}[2]{\left\langle{#1},{#2}\right\rangle}
\renewcommand{\th}{\ensuremath{\theta}}

\newcommand{\bmat}[1]{\begin{bmatrix} #1\end{bmatrix}}

\newcommand{\R}{\mathbb{R}}

\newcommand{\myint}{\int_{a}^{b}}
\newcommand{\myinta}[1]{\int_{a}^{#1}}
\newcommand{\myintb}[1]{\int_{#1}^{b}}

\newtheorem{ex}{\textbf{Example}}
\newtheorem{dem}{\textbf{Demonstration}}

\newcommand{\fourpi}[4]{\ensuremath{\mcl{P}{\tiny\bmat{#1,& \hspace{-3mm}#2 \\ #3,& \hspace{-3mm} \left\{#4\right\} }}}}

\newcommand{\threepi}[1]{\mcl{P}_{\{#1_i\}}}

\allowdisplaybreaks

\begin{document}

	\maketitle
	\thispagestyle{empty}
	\pagestyle{empty}

	\begin{abstract}
		In this paper, we present PIETOOLS, a MATLAB toolbox for the construction and handling of Partial Integral (PI) operators. The toolbox introduces a new class of MATLAB object, \texttt{opvar}, for which standard MATLAB matrix operation syntax (e.g. +, *, ' etc.) is defined. PI operators are a generalization of bounded linear operators on infinite-dimensional spaces that form a *-subalgebra with two binary operations (addition and composition) on the space $\R\times L_2$. These operators frequently appear in analysis and control of infinite-dimensional systems such as Partial Differential Equations (PDE) and Time-delay systems (TDS). Furthermore, PIETOOLS can: declare \texttt{opvar} decision variables, add operator positivity constraints, declare an objective function, and solve the resulting optimization problem using a syntax similar to the sdpvar class in YALMIP. Use of the resulting Linear Operator Inequalities (LOI) are demonstrated on several examples, including stability analysis of a PDE, bounding operator norms, and verifying integral inequalities. The result is that PIETOOLS, packaged with SOSTOOLS and MULTIPOLY, offers a scalable, user-friendly and computationally efficient toolbox for parsing, performing algebraic operations, setting up and solving convex optimization problems on PI operators. \end{abstract}
	\section{INTRODUCTION}\label{sec:intro}
	Linear operators on finite-dimensional spaces are defined by matrices. Linear Matrix Inequalities (LMI) provide a computational tool for analysis and control of dynamical systems in such finite dimensional spaces. Recently, the development of Partial Integral Equation (PIE) representations of PDE systems has created a framework for the extension of LMI-based methods to infinite-dimensional systems. The PIE representation encompasses a broad class of distributed parameter systems and is algebraic - eliminating the use of boundary conditions and continuity constraints~\cite{shivakumar_2019CDC,PEET2019132,das_2019CDC}. Such PIE representations have the form
	\begin{align*}
		\mcl H \dot{\mbf x}(t)+\mcl B_{d1}\dot w(t)+\mcl B_{d2}\dot u(t)&=\mcl A\mbf x(t)+\mcl B_1w(t)+\mcl B_2u(t)\\
		\hspace{-1cm}z(t)=\mcl C_1\mbf x(t)&+\mcl D_{11}w(t)+\mcl D_{12}u(t),\\
		\hspace{-1cm}y(t)=\mcl C_2\mbf x(t)&+\mcl D_{21}w(t)+\mcl D_{22}u(t)
	\end{align*}
	where the $\mcl H, \mcl A, \mcl B_{i}, \mcl C_i, \mcl D_{ij}$ are Partial Integral (PI) operators and have the form
	\begin{align*}
		&\left(\fourpi{P}{Q_1}{Q_2}{R_i}\bmat{x\\\mbf \Phi}\right)(s):= {\bmat{
				Px + \int_{-1}^{0} Q_1(s)\mbf \Phi(s)ds\\
				Q_2(s)x +\left(\mcl P_{\{R_i\}}\mbf \Phi\right)(s)
		}}.\vspace{-2mm}
	\end{align*}
	where
	\begin{align*}
		&\left(\mcl{P}_{\{R_i\}}\mbf \Phi\right)(s):= \\
		&R_0(s) \mbf \Phi(s) +\int_{-1}^s R_1(s,\theta)\mbf \Phi(\theta)d \theta+\int_s^0R_2(s,\theta)\mbf \Phi(\theta)d \theta
	\end{align*}
	
	PI operators, which also appear in partial-integro differential equations \cite{smyshlyaev2004closed}, have been studied in the past \cite{bergman2013integral,appell2000partial,gohberg2013classes}, extensively. PI operators are integral operators on the joint space of finite dimensional vectors and square integrable functions. Similar to matrices, PI operators are closed under the algebraic operations of addition, concatenation, composition and adjoint. As a result, LMI developed for analysis and control of finite-dimensional systems can be generalized to LOI defined by variables of the \texttt{opvar} class. For example, consider the LMI for optimal observer synthesis of singular systems: find $P\succ0$ and $Z$ such that 
	\begin{align*}
		\bmat{-\gamma I&- D_{11}^{\top}&-({P}{B}+ZD_{21})^{\top} H\\
			(\cdot)^{\top}&-\gamma I& C_1\\
			(\cdot)^{\top}&(\cdot)^{\top}&(\cdot)^{\top}+ H^{\top}({P}{A}+ZC_2)}				\prec0
	\end{align*}
	where the superscript $^{\top}$ stands for matrix transpose. This LMI can be generalized to an LOI~\cite{das_2019CDC}: Find $\mcl P=\fourpi{P}{Q}{Q^T}{R_i}\succ 0$ and $\mcl Z=\fourpi{Z_1}{\emptyset}{Z_2}{\emptyset}$ such that
	\begin{align*}
		\bmat{-\gamma I&-\mcl D_{11}^*&-(\mcl{P}\mcl{B}+\mcl{ZD}_{21})^*\mcl H\\
			(\cdot)^*&-\gamma I&\mcl C_1\\
			(\cdot)^*&(\cdot)^*&(\cdot)^*+\mcl H^*(\mcl{P}\mcl{A}+\mcl{ZC}_2)}				\prec0
	\end{align*}
	where the superscript $^*$ stands for operator adjoint. The goal of PIETOOLS~\cite{toolbox:pietools} is to create a convenient parser for constructing and solving LOI of this form. To this end, PIETOOLS incorporates all elements typically used for constructing LMIs in the commonly used LMI parser YALMIP \cite{Lofberg2004}. Specifically, PIETOOLS can be used to: declare PI operators; declare PI decision variables; manipulate PI objects via addition, multiplication, adjoint, and concatenation; add inequality constraints; set an objective function; and solve an LOI. 
	
	Significantly, PIETOOLS also includes scripts for conversion of linear TDS and coupled ODE-PDE models into PIEs. Currently, scripts are also included for stability analysis, $H_\infty$-gain analysis, $H_\infty$-optimal controller synthesis, and $H_\infty$-optimal observer synthesis. These Demo files and PIETOOLS itself are distributed as a free, third party MATLAB toolbox and are available online at \textit{{https://codeocean.com/capsule/0447940}}.
	
	The paper is organized as follows. In Section \ref{sec:note}, we introduce the standard notation utilized in the paper, followed by formal definition of Partial-Integral (PI) operators in Section \ref{sec:def}  followed by a demonstration of MATLAB implementation of the toolbox in Section \ref{sec:demo}. In the appendix, we briefly discuss algebraic operations related to PI operators which allow us to solve operator valued tests.
	\section{Notation}\label{sec:note}
	$\mathbb{S}^m\subset \R^{m\times m}$ is the set symmetric matrices. 
	For a normed space $X$, define $L_2^n[X]$ as the Hilbert space of square integrable $\R^n$-valued functions on $X$ with inner product $\langle x,y \rangle_{L_2} = \myint x(s)^{\top} y(s) ds$. The Sobolov spaces are denoted $W^{q,n}[X]:=\{x\in L_2^n[X] \mid \frac{\partial^k x}{\partial s^k}\in L_2^n[X] \text{ for all }k\le q \}$ with the standard Sobolov inner products.  $\mcl{A}^*$ stands for adjoint of a linear operator $\mcl{A}:L_2[X]\to L_2[X]$ with respect to standard $L_2$-inner product.
	For a given inner product space, $Z$, the operator $\mcl{P}:Z\to Z$ is positive semidefinite (denoted $\mcl{P}\succcurlyeq 0$) if $\ip{z}{\mcl{P}z}_Z\ge 0$ for all $z\in Z$. Furthermore, we say $\mcl{P}:Z\to Z$ is coercive if there exists some $\epsilon>0$ such that $\ip{z}{\mcl{P} z}_Z \geq \epsilon \Vert z \Vert_Z^2$ for all $z \in Z$. The partial derivative $\frac{\partial}{\partial s}\mbf x$ is denoted as $\mbf x_s$. Identity matrix of dimension $n\times n$ is denoted by $I_n$.
	\section{PI Operators And PI-operator valued optimization problems}\label{sec:def}
	Linear operators mapping between finite-dimensional spaces can be parametrized using matrices. Partial Integral operators (here onwards referred to as PI operators) are a generalization of a linear mapping between infinite-dimensional spaces, specifically a map from $\R^m\times L_2^n \to \R^p\times L_2^q$. These operators are frequently encountered in analysis and control of PDEs or TDSs.
	
	We define two class of PI-operators, 3-PI and 4-PI, where 3-PI operators are a special case of 4-PI operators. As the nomenclature insinuates, 3-PI operators, denoted as $\threepi{N}:L_2^m[a,b]\to L_2^n[a,b]$, are parameterized by $3$ matrix-valued functions $N_0: [a,b]\to\R^{n\times m}$ and $N_1, N_2: [a,b]\times[a,b] \to \R^{m\times n}$ which is a bounded linear operator between two normed spaces $L_2^m[a,b]$ and $L_2^n[a,b]$ endowed with standard $L_2$-inner product.
	{\small
	\begin{align}\label{eq:3pi}
		&\Big(\threepi{N}\mbf{y}\Big)(s) = \notag\\
		&N_0(s)\mbf{y}(s)+ \int_a^s N_1(s,\th)\mbf{y}(\th)d\th + \int_s^b N_2(s,\th)\mbf{y}(\th)d\th
	\end{align}}
	
	Similarly, 4-PI operators, parameterized by 4 components, are bounded linear operators between $\R^m\times L_2^n[a,b]$ and $\R^p\times L_2^q[a,b]$.
	\begin{align}\label{eq:4pi}
		&\fourpi{P}{Q_1}{Q_2}{R}\bmat{x\\\mbf{y}}(s) = \bmat{Px + \int_{a}^{b}Q_1(s)\mbf{y}(s)ds\\Q_2(s)x+\threepi{R}\mbf{y} (s)}
	\end{align}
	where $P:\R^m\to\R^p$, $Q_1:[a,b]\to\R^{p\times n}$, $Q_2:[a,b]\to\R^{q\times m}$ and $\threepi{R}:L_2^n[a,b]\to L_2^q[a,b]$.
	
	These operators frequently appear in control-related applications for linear TDS or coupled ODE-PDE systems. Linear TDS or coupled ODE-PDE systems with boundary conditions can be rewritten using PI operators (see \cite{shivakumar_2019CDC}). Stability test of such a system gives rise to an operator-valued feasibility test, as shown below.
	\begin{ex}[Feasibility]\label{ex:fea}
		\hspace{2mm}\\
		Test for the stability of a coupled ODE-PDE system, whose dynamics are governed by the equation, in PI format,
		\begin{align}\label{eq:ODE-PDE}
			\mcl{H} \dot{\mbf{x}} = \mcl{A} \mbf{x}
		\end{align}
		can be posed as an operator-valued feasibility test, shown below.
		\begin{align*}
			&\text{Find}, ~\mcl{P} \succ0, ~s.t.\\
			&\mcl{A}^*\mcl{PH}+\mcl{H}^*\mcl{PA} \preccurlyeq 0
		\end{align*}
		
		If there exists a self-adjoint coercive PI operator $\mcl{P}$, which satisfies the given constraints, then the system governed by Eq. \eqref{eq:ODE-PDE}, is stable.
	\end{ex}
	
	Another application of interest is finding the $\hinf$-norm of a coupled ODE-PDE system. This can be posed as an optimization problem minimizing the $L_2$-gain bound from inputs to outputs.
	
	\begin{ex}[Optimization]\label{ex:opt}
		\hspace{2mm}\\
		Finding $\hinf$-norm, $\gamma$, of a coupled ODE-PDE system whose dynamics are governed by the equation in PI format shown below
		\begin{align}\label{eq:ODE-PDE-hinf}
			\mcl{H} \dot{\mbf{x}} &= \mcl{A} \mbf{x} + \mcl{B}u,\notag\\
			y &= \mcl{C}\mbf{x}+\mcl{D}u,
		\end{align}
		can be posed as the following optimization problem.
		\begin{align*}
			&\text{minimize} ~\gamma, ~s.t.\\
			&\mcl{P} \succ0,\\
			&\bmat{-\gamma I&\mcl{D}^*&\mcl{B}^*\mcl{P}\mcl{H}\\\mcl{D}&-\gamma I&\mcl{C}\\ \mcl{H}^*\mcl{P}\mcl{B}&\mcl{C}^*&\mcl{A}^*\mcl{P}\mcl{H}+\mcl{H}^*\mcl{P}\mcl{A}}\preccurlyeq 0
		\end{align*}
	\end{ex}
	\vspace{2mm}
	Although the examples provided here are control-oriented, PIETOOLS is capable of solving other operator-valued feasibility or convex optimization problems, as described in Section~\ref{sec:demo}.
	

	
	
	\section{Declaration and Manipulation of \texttt{opvar} Objects}\label{sec:opvar}
	
	
	
	PIETOOLS introduces the structured \texttt{opvar} class of MATLAB object, each element of which consists of a PI operator $\fourpi{P}{Q_1}{Q_2}{R_i}:\R^m\times L_2^n[a,b]\to \R^p\times L_2^q[a,b]$. The structural elements of an \texttt{opvar} object are listed in Table \ref{tab:parameters}. The elements $P$, $Q1$, $Q2$, $R.R0$, $R.R1$, and $R.R2$ are themselves of the \texttt{pvar} class of polynomial introduced in the MULTIPOLY toolbox. Note that the MULTIPOLY toolbox is included in PIETOOLS toolbox, along with a modified version of SOSTOOLS. For this reason, the PIETOOLS path should take precedence over any path containing a preexisting version of MULTIPOLY or SOSTOOLS. In the solvers distributed with PIETOOLS, this is ensured by executing the Matlab command
	\texttt{addpath(genpath('.'))} from a file within the PIETOOLS directory.
	
	\begin{table}[!ht]
		\caption{List of properties in \texttt{opvar} class and their description}
		\begin{tabular}{c p{6.25cm}}
			\textbf{Property} & \textbf{Value}\\ [0.5ex]
			\hline\hline
			var1, var2 & A \texttt{pvar} object \\
			\hline
			P & A matrix with dimensions $p\times m$\\
			\hline
			Q1 & A matrix-valued \texttt{pvar} object in var1 with dimensions $p\times n$\\
			\hline
			Q2 & A matrix-valued \texttt{pvar} object in var1 with dimensions $q\times m$\\
			\hline
			R.R0 & A matrix-valued \texttt{pvar} object in var1 with dimensions $q\times n$\\
			\hline
			R.R1, R.R2 & A matrix-valued \texttt{pvar} object in var1 and var2 with dimensions $q\times n$\\
			\hline
			I & A vector with entries \texttt{[a,b]}\\
			\hline
			dim & A matrix with values \texttt{[p,m;q,n]}\\
		\end{tabular}	
		\label{tab:parameters}
	\end{table}
	
	\texttt{opvar} variables can be defined in MATLAB in two ways. The first method is directly using the \texttt{opvar} command, which is used to define opvar objects with known properties. The other method is by declaring an \texttt{opvar} decision variable - as described in Section~\ref{sec:vars}.
	
	The command \texttt{opvar} takes in string inputs and initializes them as symbolic opvar objects with default properties. These properties can be modified using standard MATLAB assignment. The following code snippet demonstrates a simple example.
	\begin{flalign*}
		&\texttt{>> opvar P1 P2;}&\\
		&\texttt{>> P1.I = [0 1];}&\\
		&\texttt{>> P1.P = rand(2,2); P1.Q1 = rand(2,1);}&
	\end{flalign*}
	
	The above code snippet would create two \texttt{opvar} variables \texttt{P1} and \texttt{P2} with default values. Next, the interval of \texttt{P1} is changed to $[0,1]$. Finally, components $P$ and $Q_1$ are reassigned with random matrices of stated dimensions. This makes \texttt{P1} a PI operator mapping $\R^2\times L_2[a,b]\to \R^2$.
	
	In addition to defining a new class, PIETOOLS overloads MATLAB operators such as \texttt{+, *} and \texttt{'} to simplify manipulation of PI-operators. Set of PI operator is a *-subalgebra with binary operations (addition and composition). This allows operations such as addition, composition, concatenation, and adjoint to be performed in a manner similar to matrices. All these operations result in another PI operator and can be performed in MATLAB.
	\subsection{Addition Of 4-PI Operators}
	Two 4-PI operators $\fourpi{A}{B_1}{B_2}{C_i}$ and $\fourpi{L}{M_1}{M_2}{N_i}$ can be added, if they have same dimensions, to obtain another 4-PI operator $\fourpi{P}{Q_1}{Q_2}{R_i}$ where
	\begin{align*}
		P&=A+L, Q_i=B_i+M_i, R_i=C_i+N_i.
	\end{align*}
	In MATLAB, two \texttt{opvar} class objects \texttt{P1} and \texttt{P2} can be added by using the MATLAB operator \texttt{+} as shown below.
	\begin{flalign*}
		&\texttt{>> P1+P2}&
	\end{flalign*}
	
	\subsection{Composing Two 4-PI Operators}
	Two 4-PI operators $\fourpi{A}{B_1}{B_2}{C_i}$ and $\fourpi{P}{Q}{Q_2}{R_i}$ can be composed if their inner dimensions match. The composition results in a new 4-PI operator $\fourpi{\hat{P}}{\hat{Q}_1}{\hat{Q}_2}{\hat{R}_i}$, where
	{\small
	\begin{align*}
		&\hat{P} = AP + \int_a^b B_1(s)Q_2(s)\text{d}s,\quad\hat{R}_0(s) = C_0(s)R_0(s),\\
		&\hat{Q}_1(s) = AQ_1(s) + B_1(s)R_0(s)+\myintb{s}B_1(\eta)R_1(\eta,s)\text{d}\eta\\
		&\qquad+\int_a^s B_1(\eta)R_2(\eta,s)\text{d}\eta,\\
		&\hat{Q}_2(s) = B_2(s)P + C_0(s)Q_2(s) + \myinta{s}C_1(s,\eta)Q_2(\eta)\text{d}\eta\\
		&\qquad+\myintb{s}C_2(s,\eta)Q_2(\eta)\text{d}\eta,\\
		&\hat{R}_1(s,\eta) =B_2(s)Q_1(\eta)+C_0(s)R_1(s,\eta)+C_1(s,\eta)R_0(\eta)\\
		&\qquad+\myinta{\eta} C_1(s,\theta)R_2(\theta,\eta)\text{d}\theta+\int_{\eta}^{s}C_1(s,\theta)R_1(\theta,\eta)\text{d}\theta\\
		&\qquad+\myintb{s}C_2(s,\theta)R_1(\theta,\eta)\text{d}\theta,\\
		&\hat{R}_2(s,\eta) =B_2(s)Q_1(\eta)+C_0(s)R_2(s,\eta)+C_2(s,\eta)R_0(\eta)\\
		&\qquad+\myinta{s} C_1(s,\theta)R_2(\theta,\eta)\text{d}\theta+\int_{s}^{\eta}C_2(s,\theta)R_2(\theta,\eta)d\theta\\
		&\qquad+\myintb{\eta}C_2(s,\theta)R_1(\theta,\eta)\text{d}\theta.
	\end{align*}}
	In MATLAB, the composition of two \texttt{opvar} objects \texttt{P1} and \texttt{P2} is performed by using \texttt{*}.
	\begin{flalign*}
		&\texttt{>> P1*P2}&
	\end{flalign*}
	
	\subsection{Adjoint Of A 4-PI Operator}
	A 4-PI operator $\fourpi{\hat{P}}{\hat{Q}_1}{\hat{Q}_2}{\hat{R}_i}$ is the adjoint of the 4-PI operator $\fourpi{P}{Q_1}{Q_2}{R_i}$, with respect to $\R\times L_2$-inner product, if 
	\begin{align*}
		&\hat{P} = P^{\top}, &\hat{R}_0(s) = R_0^{\top}(s), \nonumber\\
		&\hat{Q}_1(s) = Q_2^{\top}(s), \nonumber &\hat{R}_1(s,\eta) = R_2^{\top}(\eta,s), \nonumber\\
		&\hat{Q}_2(s) = Q_1^{\top}(s),  &\hat{R}_2(s,\eta) = R_1^{\top}(\eta,s).
	\end{align*}
	The adjoint of an \texttt{opvar} class object \texttt{P1} can be computed using the following MATLAB syntax.
	\begin{flalign*}
		&\texttt{>> P1'}&
	\end{flalign*}
	
	\subsection{Concatenation Of 4-PI Operators}
	Horizontal and vertical concatenation of \texttt{opvar} class objects, \texttt{P1} and \texttt{P2}, with compatible dimensions can be done using the following two commands, respectively.
	\begin{flalign*}
		&\texttt{>> [P1 P2]}&\\
		&\texttt{>> [P1; P2]}&
	\end{flalign*}

	\section{Declaring \texttt{opvar} Decision Variables}\label{sec:vars}
	
	Predefined \texttt{opvar} objects can be input using the syntax as described in Section~\ref{sec:opvar}. In addition, PIETOOLS can be used to set up and solve optimization problems with \texttt{opvar} decision variables. Before declaring \texttt{opvar} variables, the optimization problem structure must be initialized. This process is inherited from the SOSTOOLS toolbox and consists of the following syntax.\\
	
	\noindent\texttt{>>T = sosprogram([s,th],gam);}\\
	
	Here \texttt{s}, \texttt{th}, and \texttt{gam} are \texttt{pvar} objects. The structured object \texttt{T} carries an accumulated list of variables and constraints and must be passed whenever an additional variable or constraint is declared. The commands \texttt{sos\_opvar} and \texttt{sos\_posopvar} both declare \texttt{opvar} objects with unknown parameters. The latter function adds the constraint that the associated PI operator be positive. The syntax for both functions are listed as follows. 
	
	\subsection{\texttt{sos\_opvar}}
	\vspace{-5mm}
	\begin{flalign*}
		&\texttt{>> [T,P] = sos\_opvar(T,dim,I,s,th,deg);}&
	\end{flalign*}
	Indefinite \texttt{opvar} decision variables can be defined using the \texttt{sos\_opvar} command. This function has six required inputs: 
	\begin{enumerate}
		\item An empty or partially complete problem structure \texttt{T} to which to add the variable; 
		\item A length two vector \texttt{I=[a,b]}, indicating the spatial domain of the operator;
		\item Two \texttt{pvar} objects \texttt{s} and \texttt{th}, corresponding to the \texttt{pvar} objects declared in \texttt{sosprogram} when \texttt{T} was initialized;
		\item A $2\times 2$ matrix \texttt{dim=[p m; q n]}, indicating the dimension of domain and range of the operator; Note that when $q=n=1$, the decision variable is a matrix.
		\item A length 3 vector \texttt{deg=[d1,d2,d3]} which control the degrees of the polynomials in \texttt{P}. The value \texttt{d1} is the highest degree of \texttt{s} in the polynomials \texttt{P.Q1}, \texttt{P.Q2} and \texttt{P.R.R0}. The values \texttt{d2} and \texttt{d3} correspond to the highest degree of \texttt{s} and \texttt{th}, respectively, in the polynomials \texttt{P.R.R1} and \texttt{P.R.R2}.
	\end{enumerate}
	\texttt{sos\_opvar} returns an \texttt{opvar} object \texttt{P} which is the 4-PI operator $\fourpi{\texttt{P.P}}{\texttt{P.Q1}}{\texttt{P.Q2}}{\texttt{P.R.Ri}}:\R^m\times L_2^n[a,b] \to \R^p\times L_2^q[a,b]$ and the problem structure \texttt{T} to which the variable has been appended. For more details on construction of an \texttt{opvar} object, see Appendix E in \cite{shivakumar2019pietools}.
	\subsection{\texttt{sos\_posopvar}}
	\vspace{-5mm}
	\begin{flalign*}
		&\texttt{>> [T,P] = sos\_posopvar(T,dim,I,s,th,deg);}&
	\end{flalign*}
	Positive semi-definite \texttt{opvar} decision variables can be defined using the \texttt{sos\_posopvar} command. This function has six required inputs: 
	\begin{enumerate}
		\item An empty or partially complete problem structure \texttt{T} to which to add the variable; 
		\item A length two vector \texttt{I=[a,b]}, indicating the spatial domain of the operator;
		\item Two \texttt{pvar} objects \texttt{s} and \texttt{th}, corresponding to the \texttt{pvar} objects declared in \texttt{sosprogram} when \texttt{T} was initialized;
		\item A $2\times 1$ vector \texttt{dim=[m; n]}, indicating the dimension of domain and range of the operator. Note that when $n=0$, this becomes a standard positive matrix variable.
		\item An array \texttt{deg=[d1,d2,d3]} which control the degrees of the polynomials in \texttt{P}. The highest degree in the polynomials \texttt{P.Qi} are all max$(\texttt{d1},\texttt{d2+d3+1})$. The highest degree in \texttt{P.R.R0} is \texttt{2d2}. The highest degrees of \texttt{s} and \texttt{th} in \texttt{P.R.R1} and \texttt{P.R.R2} are all max$(\texttt{d1+d2},\texttt{2d2+d3+1})$.
	\end{enumerate}
	\texttt{sos\_posopvar} returns an \texttt{opvar} object \texttt{P} which is a self-adjoint 4-PI operator $\fourpi{\texttt{P.P}}{\texttt{P.Q1}}{\texttt{P.Q2}}{\texttt{P.R.Ri}}$ on $\R^m\times L_2^n[a,b]$ and the problem structure \texttt{T} to which the variable has been appended. The functions \texttt{P.P}, \texttt{P.Q1}, \texttt{P.Q2}, \texttt{P.R.Ri} are constrained such that \texttt{P} is a positive semidefinite operator. For more details on construction of a positive \texttt{opvar} object, see Appendix D in \cite{shivakumar2019pietools}.
	
	\section{Constraining \texttt{opvar} Class Objects And Solving An Optimization Problem}
	In addition to declaring \texttt{opvar} objects with unknown variables by using \texttt{sos\_opvar} or \texttt{sos\_posopvar}, the user can add equality constraints or operator positivity constraints to a problem structure. 
	
	\subsection{\texttt{sos\_opineq}}
	\vspace{-5mm}
	\begin{flalign*}
		&\texttt{>> T = sos\_opineq(T,P);}&
	\end{flalign*}
	\texttt{sos\_opineq} adds an operator inequality constraint of the form $\fourpi{\texttt{P.P}}{\texttt{P.Q1}}{\texttt{P.Q2}}{\texttt{P.R.Ri}}\ge 0$ to a problem structure \texttt{T}. The function has two required inputs: a problem structure \texttt{T} to which to append the constraint; and an \texttt{opvar} structure \texttt{P} which is constrained to be positive semidefinite in the augmented problem structure \texttt{T} returned by the function.
	
	\subsection{\texttt{sos\_opeq}}
	\vspace{-5mm}
	\begin{flalign*}
		&\texttt{>> T = sos\_opeq(T,P);}&
	\end{flalign*}
	\texttt{sos\_opeq} 
adds an operator equality constraint of the form $\fourpi{\texttt{P.P}}{\texttt{P.Q1}}{\texttt{P.Q2}}{\texttt{P.R.Ri}}= 0$ to a problem structure \texttt{T}. The function has two required inputs: a problem structure \texttt{T} to which to append the constraint; and an \texttt{opvar} structure \texttt{P} which is constrained to be zero in the augmented problem structure \texttt{T} returned by the function.

	For example, the constraint \texttt{P2=P1} can be imposed by using the command
	\begin{flalign*}
		&\texttt{>> T = sos\_opeq(T,P2-P1);}&
	\end{flalign*}
	where \texttt{P1} and \texttt{P2} are \texttt{opvar} objects.
	
\subsection{Defining the objective via \texttt{sossetobj}}
    Solving optimization problems using PIETOOLS, as introduced in Section \ref{sec:vars}, may require specification of an objective. This is done by using \texttt{sossetobj} function inherited from SOSTOOLS.
    \begin{flalign*}
        &\texttt{>>T = sossetobj(T,gam);}&
    \end{flalign*}
 	\texttt{sossetobj} adds the objective function \texttt{gam} to the problem structure \texttt{T}. There are two necessary inputs: The scalar \texttt{pvar} object \texttt{gam} object which is to be the minimized and the problem structure \texttt{T} to which the objective is to be added .
    
\subsection{Solving the optimization problem}
    Once all elements of the optimization problem have been added to the problem structure \texttt{T}, the problem can be solved by using the function \texttt{sossolve}, inherited from SOSTOOLS and which requires an instance of SeDuMi available in the Matlab path.
    \begin{flalign*}
    &\texttt{>>T = sossolve(T);}&
    \end{flalign*}
\texttt{sossolve} has a single input which is an `unsolved' problem structure \texttt{T}. The function returns the problem structure in a `solved' state. Data on the solution can now be obtained from the problem structure. For more details on \texttt{sossetobj} and \texttt{sossolve} we refer to the most recent SOSTOOLS documentation in~\cite{papachristodoulou2013sostools}.
    
    \subsection{\texttt{sosgetsol\_opvar}}

After execution of \texttt{sossolve}, optimal values of the real-valued decision variables may be extracted from the problem structure using the command \texttt{sosgetsol}, as described in the SOSTOOLS documentation. To extract a feasible \texttt{opvar} decision variable, the \texttt{sosgetsol\_opvar} function may be used. 

    \begin{flalign*}
    &\texttt{>>P = sosgetsol\_opvar(T,P);}&
	\end{flalign*}

\texttt{sosgetsol\_opvar} This function takes necessary inputs: a solved optimization problem structure \texttt{T} and the name of the \texttt{opvar} decision variable \texttt{P} whose solution is to be retrieved. The function returns an \texttt{opvar} object with no decision variables. This object may be manipulated further or used in the definition of a new problem structure. problem structures in the `solved' state cannot be re-used. 
	\begin{table}[!ht]
		\caption{List of functions for \texttt{opvar} class and their description}
		\begin{tabular}{l p{6.cm}}
			\textbf{Function} &\textbf{Description}\\ [0.5ex]
			\hline\hline
			opvar & Creates default opvar object with given names\\
			\hline
			+ & Adds two opvar objects\\
			\hline
			* & Composes two opvar objects\\
			\hline
			' & Transposes an opvar object\\
			\hline
			sos\_opvar & Returns opvar variable of given dimensions\\
			\hline
			sos\_posopvar & Returns a self-adjoint, opvar variable which is constrained to be positive semidefinite\\
			\hline
			sos\_opeq & Takes an input opvar variable, P, and adds the constraint P = 0 \\
			\hline
			sos\_opineq & Takes an input opvar variable, P, and adds the constraint P $\ge$ 0 \\
			\hline
			sosgetsol\_opvar & Returns the value of an opvar decision variable after solving the optimization problem
		\end{tabular}	
		\label{tab:functions}
	\end{table}

	\section{PIETOOLS scripts for analysis and control PDEs and Systems with Delay}
	As described in the Section \ref{sec:intro}, PDEs and Delay Systems admit PIE representations which can be used to test for stability, find the $\hinf$-norm or design the $\hinf$-optimal observers and controllers - See~\cite{shivakumar_2019CDC}. PIETOOLS includes the scripts \texttt{solver\_PIETOOLS\_PDE} and \texttt{solver\_PIETOOLS\_TDS}, which take input parameters as described in the header of the file and constructs to the corresponding PIE representation using the script \texttt{setup\_PIETOOLS\_PDE} or \texttt{setup\_PIETOOLS\_TDS}. Once converted to PIE form, the solver file calls one of the following executives based on the user input.
    \begin{enumerate}
        \item \texttt{executive\_PIETOOLS\_stability}: This executive is called if the user sets \texttt{stability=1} in the solver file. This tests if the PDE or TDS in PIE form is stable.
        \item \texttt{executive\_PIETOOLS\_Hinf\_gain}: This executive is called if the user sets \texttt{Hinf\_gain=1} in the solver file. The executive returns a bound on the $\hinf$-gain of the PDE or TDS in PIE form.
        \item \texttt{executive\_PIETOOLS\_Hinf\_estimator}: This executive is called if the user sets \texttt{Hinf\_estimator=1} in the solver file. The executive searches for an $\hinf$-optimal observer for the PDE or TDS in PIE form.
        \item \texttt{executive\_PIETOOLS\_Hinf\_controller}: This executive is called if the user sets \texttt{Hinf\_control=1} in the solver file. The executive searches for an $\hinf$-optimal controller for the PDE or TDS in PIE form.
    \end{enumerate}
    
For example, consider the stability test for a linear PDE system using PIETOOLS.
	\begin{dem}[Stability]
		Solution, $u$, of a tip-damped wave equation is governed by
		\begin{align*}
			u_{tt}(s,t) = u_{ss}(s,t),\quad u(0,t)=0, \quad u_s(1,t)=-ku_t(1,t).
		\end{align*}
		With a simple change of variable, this can be converted to two PDEs first-order differential in time
		\begin{align*}
			&\bmat{u_1\\u_2}_t(s,t) = \bmat{0~1\\1~0}\bmat{u_1\\u_2}_s(s,t), \\
			&u_2(0,t)=0, ~u_1(1,t)=-ku_2(1,t)
		\end{align*}
		where $u_1 = u_s$ and $u_2 = u_t$.
		
		To test for stability of this system, we seclare the parameters \texttt{a,b,A1,B} as follows in the file \texttt{solver\_PIETOOLS\_PDE} and set \texttt{stability=1}.
		\begin{flalign*}
			&\texttt{a=0; b=1;}&\\
			&\texttt{A1 = [0,1;1,0]; B = [0,1,0,0;0,0,k,1]}&
		\end{flalign*}
	\end{dem}
	\section{Solving Optimization Problems Using PIETOOLS - A Summary}
	
	As discussed in section \ref{sec:def}, PIETOOLS can solve feasibility tests or optimization problems involving \texttt{opvar} decision variables and equality/inequality constraints. This section provides a brief outline of the steps necessary for setting up and solving such an optimization problem.
	
Recall from Example \ref{ex:opt} that to find $\hinf$-norm of the system defined by Eqn.~\eqref{eq:ODE-PDE-hinf} we may solve the following optimization problem.
	\begin{align*}
		&\text{minimize} ~\gamma, ~s.t.\\
		&\mcl{P} \succ0,\\
		&\bmat{-\gamma I&\mcl{D}^*&\mcl{B}^*\mcl{P}\mcl{H}\\\mcl{D}&-\gamma I&\mcl{C}\\ \mcl{H}^*\mcl{P}\mcl{B}&\mcl{C}^*&\mcl{A}^*\mcl{P}\mcl{H}+\mcl{H}^*\mcl{P}\mcl{A}}\preccurlyeq 0
	\end{align*}
	To solve this optimization problem using PIETOOLS, the following steps are necessary.
	\begin{enumerate}
		\item Define \texttt{pvar} objects.
		\begin{flalign*}
			&\texttt{pvar s,th,gam;}&
		\end{flalign*}
		\item Initialize a problem structure. 
		\begin{flalign*}
			&\texttt{T = sosprogram([s,th],gam);}&
		\end{flalign*}
		\item Define relevant \texttt{opvar} data-objects.
		\begin{flalign*}
			&\texttt{opvar A,B,C,D,H;}&\\
			&\texttt{A=..;B=..;C=..;,D=..;H=..;X=[0,1];}&
		\end{flalign*}
		\item Add decision variables to the problem structure.
		\begin{flalign*}
			&\texttt{[T,P] = sos\_posopvar(T,dim,X,s,th);}&
		\end{flalign*}
		\item Add constraints to the problem structure.
		\begin{flalign*}
			&\texttt{D = [-gam*I \hspace{1mm} D' \hspace{5.5mm} B'*P*H;}&\\
			&\texttt{\hspace{9mm} D \hspace{8mm} -gam*I \hspace{0mm} C;}&\\
			&\texttt{\hspace{9mm} H'*P*B \hspace{0mm} C' \hspace{6mm} A'*P*H+H'*P*A];}&\\
			&\texttt{T = sosopineq(T,D);}&
		\end{flalign*}		\item Add an objective to the problem structure (minimize \texttt{gam}).
		\begin{flalign*}
			&\texttt{T = sossetobj(T,gam);}&
		\end{flalign*}
		\item Solve the completed problem structure.
		\begin{flalign*}
			&\texttt{T = sossolve(T);}&
		\end{flalign*}
		\item Extract the solution to the `solved' problem.
		\begin{flalign*}
			&\texttt{P\_s = sosgetsol\_opvar(T,P);}&
		\end{flalign*}
	\end{enumerate}
	
	\section{Demonstrations of PIETOOLS Usage}\label{sec:demo}
	In this section, a few simple examples are presented to demonstrate the use of PIETOOLS. Apart from control-related applications, described in previous sections, users can set up and solve other convex optimization problems that involve \texttt{opvar} variables. For instance, one can:  find a tight upper bound on the induced norm of a PI operator  - operators which appear in e.g. the backstepping transformation~\cite{krstic2008boundary} and input-output maps of non-linear ODEs \cite{feijoo2005associated}. Such bounds on the induced norm are obtained as follows.
	\begin{dem}[Operator norm]
		Find the $L_2$ induced operator norm for Volterra integral operator
		\begin{align*}
			(\mcl{A}x)(s) = \int_{0}^{s} x(t) dt.
		\end{align*}
		
		The operator $\mcl{A}$ is a 4-PI operator with $R_1=1$ and all other elements $0$. The $L_2$ induced operator norm is defined as $min \{\sqrt{\gamma}\mid \ip{\mcl Ax}{\mcl Ax}\leq \gamma \ip{x}{x}, \forall x\in L_2[a,b]\}$.
		The corresponding optimization problem is
		\begin{align*}
			&\text{min} ~\gamma, s.t.\\
			&\mcl{A}^*\mcl{A}\leq \gamma.
		\end{align*}
		Start by defining relevant pvar and \texttt{opvar} objects.
		\begin{flalign*}
			&\texttt{>> pvar s th gam;}&\\
			&\texttt{>> opvar A; A.R.R1 = 1; A.I=[0,1]}&
		\end{flalign*}
		
		Next, initialize a problem structure with $s$, $th$ and $gam$ as \texttt{pvar} objects. $gam$ is the objective to be minimized.
		\begin{flalign*}
			&\texttt{>> prog = sosprogram([s,th],[gam]);}&\\
			&\texttt{>> prog = sossetobj(prog,gam);}&
		\end{flalign*}
		Next, add the opvar inequality constraint  using the \texttt{sos\_opineq} function.
		\begin{flalign*}
			&\texttt{>> prog = sos\_opineq(prog, A'*A-gam);}&
		\end{flalign*}
		Finally, the problem can be solved and solution extracted using the following commands.
		\begin{flalign*}
			&\texttt{>> prog = sossolve(prog);}&\\
			&\texttt{>> Gam = sosgetsol(prog, gam);}&\\
			&\texttt{>> disp(sqrt(Gam));}&\\
			&\texttt{ans =}&\\
			&\texttt{ \quad \quad 0.6366}&
		\end{flalign*}
		The numerical value of $.6366$ obtained from PIETOOLS can be compared to the analytical value of the induced norm of this operator norm which is known to be $\frac{2}{\pi}\approx 0.6366$.
	\end{dem}
	
	PIETOOLS can also be used to provide certificates of positivity for integral inequalities. 
	
	\begin{dem}[Poincare's Constant]
		Poincare's Inequality states that there exists a constant $C$ such that for every function $u\in W^{1,p}_0(\Omega)$ (where $W^{1,p}_0(\Omega)$ is the Sobolev space of zero-trace functions) we have that
		\begin{align*}
			\norm{u}_{L_p(\Omega)}\leq C \norm{\nabla u}_{L_p(\Omega)}
		\end{align*} where $1\leq p < \infty$ and $\Omega$ is a bounded set.
		This can be rewritten as an optimization problem.
		\begin{align*}
			&\text{min} ~C, s.t.\\
			&\ip{u}{u} - C\ip{u_s}{u_s}\leq 0.
		\end{align*}
		For $p=2$ and $\Omega = [0,1]$, it known that  for functions $u\in W^{2,2}_0(\Omega)$ with the boundary conditions $u(0)=u(1)=0$ the smallest  $C= 1/\pi^2$. Numerical calculation of this constant $C$ can be reformulated as a PI optimization problem as follows.
		\begin{align*}
			&\text{min} ~C, s.t.\\
			&\mcl{H^*H}-C\mcl{H}_2^*\mcl{H}_2\leq 0\\
			& (\mcl{H}u) (s) = \myinta{s} (s\th -\th) u(\th)d\th + \myintb{s} (s\th -s) u(\th)d\th\\
			&(\mcl{H}_2 u)(s) =\myinta{s} \th u(\th)d\th + \myintb{s} (\th -1) u(\th)d\th
		\end{align*}
		
		The set up and solution of this PI optimization problem using PIETOOLS is as follows.
		\begin{flalign*}
			&\texttt{pvar s t C; opvar H H2;}&\\
			&\texttt{H.I=[0,1]; H.R.R1=s*t-t; H.R.R2=s*t-s;}&\\
			&\texttt{H2.I=[0,1]; H2.R.R1 = t; H2.R.R2 = t-1;}&\\
			&\texttt{prog = sosprogram([s,t],C);}&\\
			&\texttt{prog = sossetobj(prog,C);}&\\
			&\texttt{prog = sos\_opineq(prog, H'*H-C*H2'*H2);}&\\
			&\texttt{prog = sossolve(prog);}&
		\end{flalign*}
	When implemented, this code returns a smallest bound of $\sqrt{C}=.3183\approx \frac{1}{\pi}$
	\end{dem}
	
	%
	\section{Conclusions}
	In this paper, we have provided a guide to the new MATLAB toolbox PIETOOLS for manipulation and optimization of PI operators. We have provided details on declaration of PI operator objects, manipulation of PI operators, declaration of PI decision variables, addition of operator equality and inequality constraints, solution of PI optimization problems, and extraction of feasible operators. We have demonstrated the practical usage of PIETOOLS, including scripts for analysis and control of PDEs and systems with delay, as well as bounding operator norms and proving integral inequalities. These examples and descriptions illustrate both the syntax of available features and the necessary components of any PIETOOLS script. Finally, we note that PIETOOLS is still under active development. Ongoing efforts focus on identifying and balancing the degree structures in \texttt{sos\_opineq}.

	\section*{ACKNOWLEDGMENT}
	This work was supported by Office of Naval Research Award N00014-17-1-2117 and National Science Foundation under Grants No. 1739990 and 1935453.
	
	\bibliographystyle{IEEEtran}
	\bibliography{references,peet_bib}

\begin{thebibliography}{10}
\providecommand{\url}[1]{#1}
\csname url@samestyle\endcsname
\providecommand{\newblock}{\relax}
\providecommand{\bibinfo}[2]{#2}
\providecommand{\BIBentrySTDinterwordspacing}{\spaceskip=0pt\relax}
\providecommand{\BIBentryALTinterwordstretchfactor}{4}
\providecommand{\BIBentryALTinterwordspacing}{\spaceskip=\fontdimen2\font plus
\BIBentryALTinterwordstretchfactor\fontdimen3\font minus
  \fontdimen4\font\relax}
\providecommand{\BIBforeignlanguage}[2]{{%
\expandafter\ifx\csname l@#1\endcsname\relax
\typeout{** WARNING: IEEEtran.bst: No hyphenation pattern has been}%
\typeout{** loaded for the language `#1'. Using the pattern for}%
\typeout{** the default language instead.}%
\else
\language=\csname l@#1\endcsname
\fi
#2}}
\providecommand{\BIBdecl}{\relax}
\BIBdecl

\bibitem{shivakumar_2019CDC}
S.~Shivakumar, A.~Das, S.~Weiland, and M.~Peet, ``A generalized {LMI}
  formulation for input-output analysis of linear systems of {ODE}s coupled
  with {PDEs},'' in \emph{Proceedings of the IEEE Conference on Decision and
  Control}, 2019.

\bibitem{PEET2019132}
M.~M. Peet, S.~Shivakumar, A.~Das, and S.~Weiland, ``Discussion paper: A new
  mathematical framework for representation and analysis of coupled {PDEs},''
  \emph{3rd IFAC Workshop on Control of Systems Governed by Partial
  Differential Equations CPDE 2019}, vol.~52, no.~2, pp. 132 -- 137, 2019.

\bibitem{das_2019CDC}
A.~Das, S.~Shivakumar, S.~Weiland, and M.~Peet, ``${H}_\infty$ optimal
  estimation for linear coupled {PDE} systems,'' in \emph{Proceedings of the
  IEEE Conference on Decision and Control}, 2019.

\bibitem{smyshlyaev2004closed}
A.~Smyshlyaev and M.~Krstic, ``Closed-form boundary state feedbacks for a class
  of {1D} partial integro-differential equations,'' \emph{IEEE Transactions on
  Automatic Control}, vol.~49, no.~12, pp. 2185--2202, 2004.

\bibitem{bergman2013integral}
S.~Bergman, \emph{Integral operators in the theory of linear partial
  differential equations}.\hskip 1em plus 0.5em minus 0.4em\relax Springer,
  2013.

\bibitem{appell2000partial}
J.~Appell, A.~Kalitvin, and P.~Zabrejko, \emph{Partial Integral Operators and
  Integro-Differential Equations: Pure and Applied Mathematics}.\hskip 1em plus
  0.5em minus 0.4em\relax CRC Press, 2000.

\bibitem{gohberg2013classes}
I.~Gohberg, S.~Goldberg, and M.~A. Kaashoek, \emph{Classes of linear
  operators}.\hskip 1em plus 0.5em minus 0.4em\relax Birkh{\"a}user, 2013,
  vol.~63.

\bibitem{toolbox:pietools}
M.~Peet and S.~Shivakumar, ``{PIETOOLS} for {Networks} with {Delay},''
  \url{https://codeocean.com/capsule/0447940}.

\bibitem{Lofberg2004}
J.~L{\"{o}}fberg, ``{YALMIP} : A toolbox for modeling and optimization in
  {MATLAB},'' in \emph{In Proceedings of the CACSD Conference}, Taipei, Taiwan,
  2004.

\bibitem{shivakumar2019pietools}
S.~Shivakumar, A.~Das, and M.~M. Peet, ``{PIETOOLS}: A {MATLAB} toolbox for
  manipulation and optimization of partial integral operators,'' \emph{arXiv
  preprint arXiv:1910.01338}, 2019.

\bibitem{papachristodoulou2013sostools}
A.~Papachristodoulou, J.~Anderson, G.~Valmorbida, S.~Prajna, P.~Seiler, and
  P.~Parrilo, ``{SOSTOOLS} version 3.00 sum of squares optimization toolbox for
  {MATLAB},'' \emph{arXiv preprint arXiv:1310.4716}, 2013.

\bibitem{krstic2008boundary}
M.~Krstic and A.~Smyshlyaev, \emph{Boundary control of {PDEs}: A course on
  backstepping designs}.\hskip 1em plus 0.5em minus 0.4em\relax Siam, 2008,
  vol.~16.

\bibitem{feijoo2005associated}
J.~V. Feijoo, K.~Worden, and R.~Stanway, ``Associated linear equations for
  {Volterra} operators,'' \emph{Mechanical Systems and Signal Processing},
  vol.~19, no.~1, pp. 57--69, 2005.

\end{thebibliography}

\end{document}